# Finding the Sequence of Largest Small $n$-Polygons by Numerical Optimization

János D. Pintér [a], Frank J. Kampas [b], Ignacio Castillo [c, *]

[a] Department of Management Science and Information Systems, Rutgers University, Piscataway, NJ, USA
[b] Physicist at Large Consulting LLC, Bryn Mawr, PA, USA
[c] Lazaridis School of Business and Economics, Wilfrid Laurier University, Waterloo, ON, Canada
[*] Corresponding author. E-mail address: icastillo@wlu.ca

**Abstract**

$LSP(n)$, the *largest small polygon* with $n$ vertices, is the polygon of unit diameter that has maximal area $A(n)$. It is known that for all *odd* values $n \geq 3$, $LSP(n)$ is the regular $n$-polygon; however, this statement is not valid for *even* values of $n$. Finding the polygon $LSP(n)$ and $A(n)$ for even values $n \geq 6$ has been a long-standing challenge. In this work, we develop high-precision numerical solution estimates of $A(n)$ for even values $n \geq 4$, using the *Mathematica* model development environment and the IPOPT local nonlinear optimization solver engine. First, we present a revised (tightened) LSP model that greatly assists the efficient solution of the model-class considered. This is followed by numerical results for an illustrative sequence of even values of $n$, up to $n \leq 1000$. Our results are in close agreement with, or surpass, the best results reported in all earlier studies. Most of these earlier works addressed special cases up to $n \leq 20$, while others obtained numerical optimization results for a range of values from $6 \leq n \leq 100$. For completeness, we also calculate numerically optimized results for a selection of odd values of $n$, up to $n \leq 999$: these results can be compared to the corresponding theoretical (exact) values. The results obtained are used to provide regression model-based estimates of the optimal area sequence $\{A(n)\}$, for all even and odd values $n$ of interest, thereby essentially solving the entire LSP model-class *numerically*, with demonstrably high precision.

**Keywords**

Nonlinear programming; Largest small polygons (LSP); $\{LSP(n)\}$ model-class; Optimal area sequence $\{A(n)\}$; Revised LSP model; *Mathematica* model development environment; IPOPT solver engine; Numerical optimization results and regression models for estimating $\{A(n)\}$.

## 1 Introduction

The *diameter* of a (convex planar) polygon is defined as the length of its longest diagonal. The *largest small polygon* (LSP) with $n$ vertices is defined as the polygon of unit diameter that has maximal area. For a given integer $n \geq 3$, we refer to this polygon as $LSP(n)$ with corresponding area $A(n)$. To illustrate, Figure 1 shows visual representations of $LSP(6)$ and $LSP(18)$. $LSP(6)$ is cited from the *MathWorld* website of Weisstein (2020), while $LSP(18)$ is based on our result:

we will present further $LSP(n)$ configuration examples later in this work. All diagonals of unit length are shown in red, all others are shown in black.

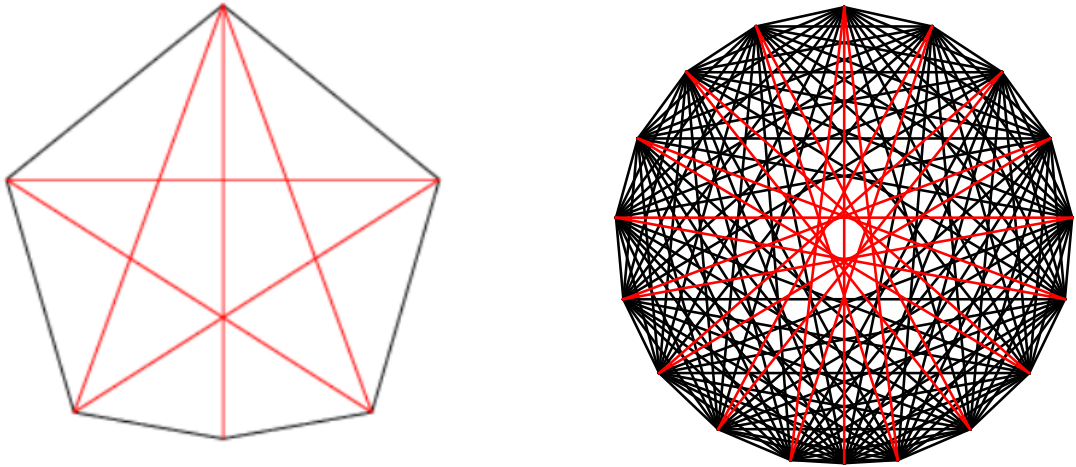

Figure 1. The largest small polygons $LSP(6)$ and $LSP(18)$.

Nearly a century ago, Reinhardt (1922) proved that for all *odd* values $n \geq 3$, $LSP(n)$ is the regular $n$-polygon; perhaps surprisingly, the corresponding statement is not valid for *even* values of $n$. For brevity, here we only refer to the discussion of the low-dimensional cases studied elsewhere. Specifically, Mossinghoff (2006) presents some related theoretical background, and confirms earlier best-known results up to $n = 16$, extending the range to $6 \leq n \leq 20$. In a recent study, Pintér (2020) gives globally optimized numerical estimates of $A(n)$ for all even values $6 \leq n \leq 80$, with comparative references to the earlier works. These works clearly indicate the limitations and varying performance level of the modeling and optimization software packages used earlier.

In this work, we present the results of an in-depth numerical optimization study aimed at finding putative $LSP(n)$ configurations and corresponding values $A(n)$, for a substantial selection of even values covering the range $4 \leq n \leq 1000$. We present a revised (tightened) LSP model, then solve model instances using the *Mathematica* model development platform with the callable IPOPT nonlinear optimization solver. Our numerical results closely match or surpass all the *best* results attained individually by several groups of researchers. For completeness, we also calculate numerically optimized results for selected odd values $3 \leq n \leq 999$ which can be compared to their theoretical values (easily estimated by computation). Based on the results obtained, we also present two regression models to estimate $\{A(n)\}$, separately for all relevant even and odd values of $n$. Since for large $n$ the actual calculated optimum estimates $A(n)$ already approximate closely the theoretical limit $A(\infty) = \pi/4$, our study enables the high-precision handling of the entire LSP model-class *numerically*.

## 2 Revised LSP Optimization Model

For unambiguity, we consider all $LSP(n)$ instances for even $n \geq 4$ (for odd $n \geq 3$) with a fixed position corresponding to appropriate variants of the instances shown by Figure 1. Following standard postulated assumptions, each even $n$-polygon considered here is symmetrical with respect to the diagonal that connects its lowest positioned vertex placed at the origin with its highest vertex. We refer to this configuration as the standard model: it has been used in all topical works referenced by Pintér (2020).

The standard LSP model uses polar coordinates to describe $LSP(n)$: vertex $i$ is positioned at polar radius $r_i$ and at angle $\theta_i$. For unambiguity, we postulate that the vertices $i = 1, \ldots, n - 1$ are arranged according to increasing angles $\theta_i$. Placing the last vertex at the origin, we set $r_n = 0$, $\theta_n = \pi$: recall Figure 1 that corresponds to such a configuration for $LSP(6)$ and $LSP(18)$. The corresponding standard LSP optimization model is presented below.

Maximize total area of the $n$-polygon:
$$\max A(n) = \frac{1}{2} \sum_{i=1}^{n-1} r_i r_{i+1} \sin(\theta_{i+1} - \theta_i). \tag{1}$$

Prescribed upper bound for the pairwise distance between vertices $i$ and $j$:
$$r_i^2 + r_j^2 - 2 r_i r_j \cos(\theta_i - \theta_j) \leq 1, \text{ for } 1 \leq i \leq n - 2, i + 1 \leq j \leq n - 1. \tag{2}$$

Vertex angle ordering relations:
$$\theta_{i+1} - \theta_i \geq 0, \text{ for } 1 \leq i \leq n - 2. \tag{3}$$

Decision variable bounds, and the two fixed variable settings:
$$0 \leq \theta_i \leq \pi \text{ and } 0 \leq r_i \leq 1, \text{ for } 1 \leq i \leq n - 1; r_n = 0, \theta_n = \pi. \tag{4}$$

Based on the structure of LSP configurations found in all earlier studies, next we revise this standard model, by adding the relations shown below.

ii) We postulate bounds on the angle differences, for even $n$:
$$\theta_{i+1} - \theta_i \geq \frac{\pi}{n}, \text{ for } 1 \leq i \leq n - 2, \tag{5}$$
$$\theta_{n/2} = \frac{\pi}{2}. \tag{6}$$

For odd $n$:
$$\theta_{i+1} - \theta_i = \frac{\pi}{n}, \text{ for } 1 \leq i \leq n - 2, \tag{7}$$
$$\theta_1 + \theta_{n-1} = \pi. \tag{8}$$

ii) We postulate the symmetry of the LSP configuration to be found, for even $n$:
$$r_{n/2+i-1} = r_{n/2-i+1}, \text{ for } 2 \leq i \leq \frac{n}{2}, \tag{9}$$
$$\theta_{n/2+i-1} = \pi - \theta_{n/2-i+1}, \text{ for } 2 \leq i \leq \frac{n}{2}. \tag{10}$$

To illustrate these added constraints, recall Figure 1; further visual examples will be presented later on. Our preliminary experimentation demonstrated that the symmetry postulates (9)-(10) for even $n$, despite reducing the number of decision variables, were not useful within our numerical optimization study. However, the bound postulates on the angle differences ((5)-(6) for even $n$, (7)-(8) for odd $n$) were highly useful, effectively tightening the LSP model.

The added bound postulates on the angle differences make the LSP model tighter, and thereby – as expected – easier to solve. As it turns out, the new constraints are essential to guarantee the performance of the *local* solver IPOPT in solving the *global* optimization problem (1)-(4), with the added considerations (5)-(6) for even $n$, (7)-(8) for odd $n$. Our results are in close agreement with or surpass all the *best* results reported earlier.

Notice the potential numerical challenge, implied by the nonconvex objective function (1) and constraints (2): the number of these constraints increases *quadratically* as a function of $n$. As noted by other researchers including Pintér (2020), while the standard $LSP(n)$ model instances have a unique globally optimal solution, the number of local optima increases with $n$. Many of the local optima are close in quality to the (only approximately known, since numerically estimated) global optimum. These features make the $\{LSP(n)\}$ problem-class hard to solve, similarly to many other object configuration design and packing problems arising e.g. in computational mathematics, physics, chemistry, biology, as well as across a range of engineering applications.

## 3 Numerical Results for Even and Odd Values $3 \leq n \leq 1000$

To formulate the LSP model versions, we use the *Mathematica* model development environment, and the IPOPT local nonlinear optimization solver engine linked to *Mathematica*. Result analysis and visualization is also done in *Mathematica.*

*Mathematica*, by Wolfram Research, is an integrated software platform that enables the analysis of an impressive range of scientific and technical computing topics. *Mathematica* is available on all modern desktop systems, as well as in the cloud through web browsers; and it is used by millions of educators, researchers, and students around the world.

Specifically, *Mathematica* can be used to develop concise and fully scalable optimization models. For a related exposition, consult e.g. Pintér and Kampas (2005); for a more recent benchmarking study using our *MathOptimizer Professional* global-local optimization software linked to *Mathematica*, see Pintér and Kampas (2013).

In our present study, we use the IPOPT (Interior Point OPTimizer) nonlinear optimization solver package (Wächter and Laird, 2020). IPOPT can be used to find local solutions to large-scale optimization models. The *Mathematica*-IPOPT link and its usage options are documented at Wolfram Research (2020b).

Since IPOPT is a local scope solver, it requires an initial solution "guess", and it greatly benefits from a good choice of that solution. Considering the postulated structure of the LSP configurations to be found, for a given $n$ our initial angle settings are chosen as $\theta_i = i(\pi/n)$, for $1 \leq i \leq n$; together with the initial radius settings $r_i = 1$ for $1 \leq i \leq n - 1$; $r_n = 0$.

The study summarized here was conducted on a laptop PC with the following specifications: Intel Core i7-7700 CPU @ 3.6 GHz (x-64 processor), 16.0 GB RAM, running under the Windows 10 Pro (64-bit) operating system.

The numerical results obtained by using *Mathematica* and IPOPT are summarized in the tables presented below. Detailed $LSP(n)$ configurations (listing all decision variable values found and all constraints) can be optionally reported. All optimization results, with related analysis and visualization are directly cited from our corresponding *Mathematica* (notebook) documents: thereby the entire study (with identical results) is completely reproducible.

Table 1 summarizes our results for selected even values $4 \leq n \leq 1000$. Let us point out that since the solution times become substantial as $n$ increases, we report results only for a rather broad selection of even (and also for odd) values; however, in principle we could handle *all* instances of $LSP(n)$, provably at least up to $n \leq 1000$. For example (see Table 1), the runtime for $n = 10$ is only 0.06 seconds; for $n = 100$ it is still just 9.44 seconds; but for $n = 1000$ it becomes 5417.51 seconds. We did not conduct systematic tests to find the largest possible numerical instance that we could handle using *Mathematica* with IPOPT, noting that – obviously – all modeling systems and optimization engines have their inherent limitations, depending also on the hardware platform and other circumstances.

The legend used is self-explanatory, "Maximum violation" refers to the maximal constraint violation at the numerical optimal solution. One can verify the *linear* increase in the number of decision variables and the rapid *quadratic* increase in the number of constraints as a function of $n$ (recalling the LSP model). The model instance for $n = 1000$ has almost two thousand decision variables and nearly half a million constraints.

Table 1. *Mathematica*-IPOPT numerical results for a selection of even values of $n$

| $n$ | Decision variables | Constraints | Runtime (seconds) | Objective function | Maximum violation |
|---|---|---|---|---|---|
| 4 | 6 | 5 | 0.03 | 0.500000 | 9.9636E-09 |
| 6 | 10 | 14 | 0.03 | 0.674981 | 9.9432E-09 |
| 8 | 14 | 27 | 0.05 | 0.726868 | 9.9236E-09 |
| 10 | 18 | 44 | 0.06 | 0.749137 | 9.9046E-09 |
| 12 | 22 | 65 | 0.08 | 0.760730 | 9.8855E-09 |
| 14 | 26 | 90 | 0.11 | 0.767531 | 9.8663E-09 |
| 16 | 30 | 119 | 0.15 | 0.771861 | 9.8472E-09 |
| 18 | 34 | 152 | 0.17 | 0.774788 | 9.8296E-09 |
| 20 | 38 | 189 | 0.23 | 0.776859 | 9.8101E-09 |
| 24 | 46 | 275 | 0.32 | 0.779524 | 9.7801E-09 |
| 28 | 54 | 377 | 0.45 | 0.781111 | 9.7647E-09 |
| 32 | 62 | 495 | 0.68 | 0.782133 | 9.8456E-09 |
| 36 | 70 | 629 | 0.82 | 0.782828 | 9.6522E-09 |
| 40 | 78 | 779 | 1.05 | 0.783323 | 9.6132E-09 |
| 44 | 86 | 945 | 1.34 | 0.783687 | 9.5741E-09 |
| 48 | 94 | 1127 | 1.53 | 0.783964 | 9.5352E-09 |
| 52 | 102 | 1325 | 1.85 | 0.784178 | 9.4967E-09 |
| 56 | 110 | 1539 | 3.04 | 0.784252 | 9.6754E-09 |
| 60 | 118 | 1769 | 3.34 | 0.784408 | 9.6548E-09 |
| 70 | 138 | 2414 | 3.92 | 0.784729 | 9.3199E-09 |
| 80 | 158 | 3159 | 4.98 | 0.784886 | 9.2227E-09 |
| 90 | 178 | 4004 | 7.59 | 0.784994 | 9.1242E-09 |
| 100 | 198 | 4949 | 9.44 | 0.785072 | 9.0264E-09 |
| 110 | 218 | 5994 | 11.18 | 0.785129 | 8.9291E-09 |
| 120 | 238 | 7139 | 14.21 | 0.785172 | 8.8306E-09 |
| 130 | 258 | 8384 | 17.52 | 0.785205 | 8.7330E-09 |
| 140 | 278 | 9729 | 21.39 | 0.785232 | 8.7013E-09 |
| 150 | 298 | 11174 | 25.11 | 0.785254 | 8.8048E-09 |
| 160 | 318 | 12719 | 29.16 | 0.785271 | 8.4389E-09 |
| 180 | 358 | 16109 | 52.43 | 0.785298 | 8.2424E-09 |
| 200 | 398 | 19899 | 51.31 | 0.785317 | 8.0460E-09 |
| 220 | 438 | 24089 | 64.43 | 0.785331 | 7.8491E-09 |
| 240 | 478 | 28679 | 81.60 | 0.785342 | 7.6515E-09 |
| 260 | 518 | 33669 | 97.73 | 0.785350 | 7.6285E-09 |
| 280 | 558 | 39059 | 118.77 | 0.785357 | 7.2925E-09 |
| 300 | 598 | 44849 | 142.19 | 0.785362 | 7.0879E-09 |
| 320 | 638 | 51039 | 170.00 | 0.785367 | 6.8695E-09 |
| 340 | 678 | 57629 | 202.27 | 0.785370 | 6.6526E-09 |
| 360 | 718 | 64619 | 235.11 | 0.785373 | 6.4506E-09 |
| 380 | 758 | 72009 | 269.55 | 0.785376 | 6.2473E-09 |
| 400 | 798 | 79799 | 316.30 | 0.785378 | 6.0432E-09 |
| 420 | 838 | 87989 | 363.92 | 0.785380 | 5.8990E-09 |
| 440 | 878 | 96579 | 393.75 | 0.785381 | 5.6483E-09 |
| 460 | 918 | 105569 | 464.32 | 0.785383 | 5.4201E-09 |
| 480 | 958 | 114959 | 514.01 | 0.785384 | 5.2503E-09 |
| 500 | 998 | 124749 | 577.96 | 0.785385 | 5.0971E-09 |
| 550 | 1098 | 150974 | 739.51 | 0.785387 | 4.4905E-09 |
| 600 | 1198 | 179699 | 948.26 | 0.785389 | 4.0205E-09 |
| 650 | 1298 | 210924 | 1228.83 | 0.785391 | 3.4442E-09 |
| 700 | 1398 | 244649 | 1460.31 | 0.785392 | 2.9080E-09 |
| 750 | 1498 | 280874 | 1857.58 | 0.785392 | 2.3883E-09 |
| 800 | 1598 | 319599 | 2342.02 | 0.785393 | 1.8715E-09 |
| 850 | 1698 | 360824 | 3177.43 | 0.785394 | 1.3518E-09 |
| 900 | 1798 | 404549 | 3846.06 | 0.785394 | 8.2861E-10 |
| 950 | 1898 | 450774 | 3721.68 | 0.785395 | 3.0161E-10 |
| 1000 | 1998 | 499499 | 5417.51 | 0.785395 | 0.0000E+00 |

Although we could directly determine (and then numerically estimate) the optimum $A(n)$ for all odd values $n$, for comparison we conducted detailed tests also for this case. Table 2 summarizes our optimization results for a selection of odd values $3 \leq n \leq 1000$. Notice again the rapidly increasing model sizes and runtimes.

Table 2. *Mathematica*-IPOPT numerical results for a selection of odd values of $n$

| $n$ | Decision variables | Constraints | Runtime (seconds) | Objective function | Maximum violation |
|---|---|---|---|---|---|
| 3 | 4 | 2 | 0.01 | 0.433013 | 5.7718E-09 |
| 5 | 8 | 9 | 0.03 | 0.657164 | 9.9306E-09 |
| 7 | 12 | 20 | 0.04 | 0.719741 | 9.9122E-09 |
| 9 | 16 | 35 | 0.06 | 0.745619 | 9.8911E-09 |
| 11 | 20 | 54 | 0.08 | 0.758748 | 9.8685E-09 |
| 13 | 24 | 77 | 0.11 | 0.760920 | 9.9209E-09 |
| 15 | 28 | 104 | 0.13 | 0.771056 | 9.8227E-09 |
| 17 | 32 | 135 | 0.17 | 0.774230 | 9.8037E-09 |
| 19 | 36 | 170 | 0.21 | 0.774632 | 9.8870E-09 |
| 23 | 44 | 252 | 0.31 | 0.778297 | 9.8699E-09 |
| 27 | 52 | 350 | 0.45 | 0.780369 | 9.7408E-09 |
| 31 | 60 | 464 | 0.66 | 0.781646 | 9.6998E-09 |
| 35 | 68 | 594 | 0.82 | 0.782492 | 9.6589E-09 |
| 39 | 76 | 740 | 1.11 | 0.783081 | 9.6656E-09 |
| 43 | 84 | 902 | 1.26 | 0.783508 | 9.6322E-09 |
| 47 | 92 | 1080 | 1.71 | 0.783827 | 9.5993E-09 |
| 51 | 100 | 1274 | 1.94 | 0.784071 | 9.5659E-09 |
| 55 | 108 | 1484 | 2.44 | 0.784331 | 9.4267E-09 |
| 59 | 116 | 1710 | 3.45 | 0.784416 | 9.4981E-09 |
| 69 | 136 | 2345 | 3.83 | 0.784686 | 9.4154E-09 |
| 79 | 156 | 3080 | 5.65 | 0.784854 | 9.5434E-09 |
| 89 | 176 | 3915 | 7.04 | 0.784975 | 9.2480E-09 |
| 99 | 196 | 4850 | 9.69 | 0.785057 | 9.2965E-09 |
| 109 | 216 | 5885 | 11.02 | 0.785118 | 9.0806E-09 |
| 119 | 236 | 7020 | 14.30 | 0.785164 | 9.1430E-09 |
| 129 | 256 | 8255 | 17.16 | 0.785199 | 8.9132E-09 |
| 139 | 276 | 9590 | 19.91 | 0.785227 | 8.9939E-09 |
| 149 | 296 | 11025 | 24.07 | 0.785249 | 9.0230E-09 |
| 159 | 316 | 12560 | 30.41 | 0.785268 | 8.8522E-09 |
| 179 | 356 | 15930 | 38.75 | 0.785295 | 8.7010E-09 |
| 199 | 396 | 19700 | 49.02 | 0.785315 | 8.5547E-09 |
| 219 | 436 | 23870 | 64.81 | 0.785330 | 8.4089E-09 |
| 239 | 476 | 28440 | 78.86 | 0.785341 | 8.2614E-09 |
| 259 | 516 | 33410 | 133.35 | 0.785349 | 8.1288E-09 |
| 279 | 556 | 38780 | 130.45 | 0.785356 | 7.9684E-09 |
| 299 | 596 | 44550 | 147.40 | 0.785362 | 7.8237E-09 |
| 319 | 636 | 50720 | 174.64 | 0.785366 | 7.5393E-09 |
| 339 | 676 | 57290 | 213.06 | 0.785370 | 7.5307E-09 |
| 359 | 716 | 64260 | 240.20 | 0.785373 | 7.3824E-09 |
| 379 | 756 | 71630 | 282.64 | 0.785375 | 7.2357E-09 |
| 399 | 796 | 79400 | 306.82 | 0.785378 | 7.0888E-09 |
| 419 | 836 | 87570 | 364.07 | 0.785380 | 6.9439E-09 |
| 439 | 876 | 96140 | 407.24 | 0.785381 | 6.7956E-09 |
| 459 | 916 | 105110 | 456.32 | 0.785383 | 6.6498E-09 |
| 479 | 956 | 114480 | 514.04 | 0.785384 | 6.5050E-09 |
| 499 | 996 | 124250 | 582.43 | 0.785385 | 6.3570E-09 |
| 549 | 1096 | 150425 | 735.91 | 0.785387 | 5.9970E-09 |
| 599 | 1196 | 179100 | 933.16 | 0.785389 | 5.3487E-09 |
| 649 | 1296 | 210275 | 1240.52 | 0.785390 | 4.9833E-09 |
| 699 | 1396 | 243950 | 1498.28 | 0.785392 | 4.5878E-09 |
| 749 | 1496 | 280125 | 1879.67 | 0.785392 | 4.2713E-09 |
| 799 | 1596 | 318800 | 2271.86 | 0.785393 | 3.8937E-09 |
| 849 | 1696 | 359975 | 2876.72 | 0.785394 | 3.9129E-09 |
| 899 | 1796 | 403650 | 4057.11 | 0.785394 | 3.2825E-09 |
| 949 | 1896 | 449825 | 3794.15 | 0.785395 | 2.9633E-09 |
| 999 | 1996 | 498500 | 5784.43 | 0.785395 | 2.6393E-09 |

Although Tables 1 and 2 show optimized $A(n)$ values only with 6 digits after the decimal point, the actually reported precision in our detailed numerical tests (within *Mathematica*) is set to 10 digits after the decimal point. For example, our numerically optimized estimate for $A(1000)$ approximately equals $0.7853949284$. Such higher precision is in line with the required constraint satisfaction level, all in the order of $10^{-9}$ as shown in the tables. The 10-digit precision also supports in-depth comparisons with results obtained earlier. Specifically, our results are in close agreement with or surpass all the *best* results reported earlier, including Mossinghoff (2006) and Pintér (2020) with reported comparisons to all earlier topical works (known to us). Furthermore, our modeling and optimization approach enables solving $LSP(n)$ instances for significantly higher values of $n$ than done before.

The following figures display the solutions found, for a selection of even and odd sequences of $n$, showing all pairwise vertex connections; again, the diagonals of unit length are shown in red, all others are shown in black. Observe that the resulting configurations quickly approach the circle, as $n$ increases.

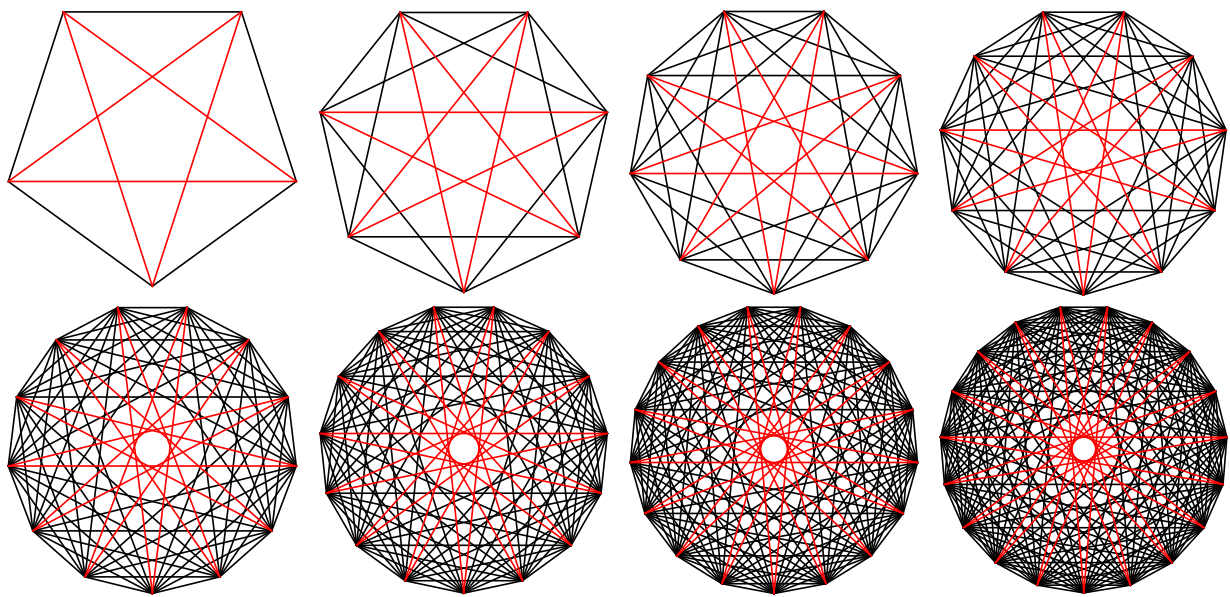

Figure 2. Largest small polygons $LSP(n)$ for $n = 5, 7, 9, 11, 13, 15, 17, 19$.

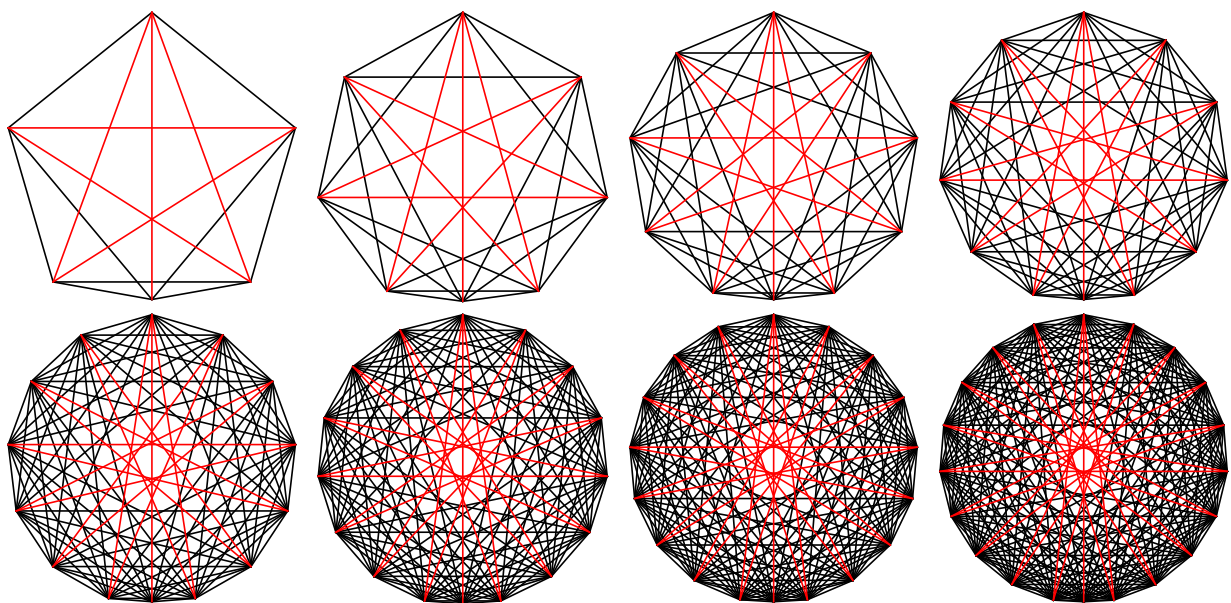

Figure 3. Largest small polygons $LSP(n)$ for $n = 6, 8, 10, 12, 14, 16, 18, 20$.

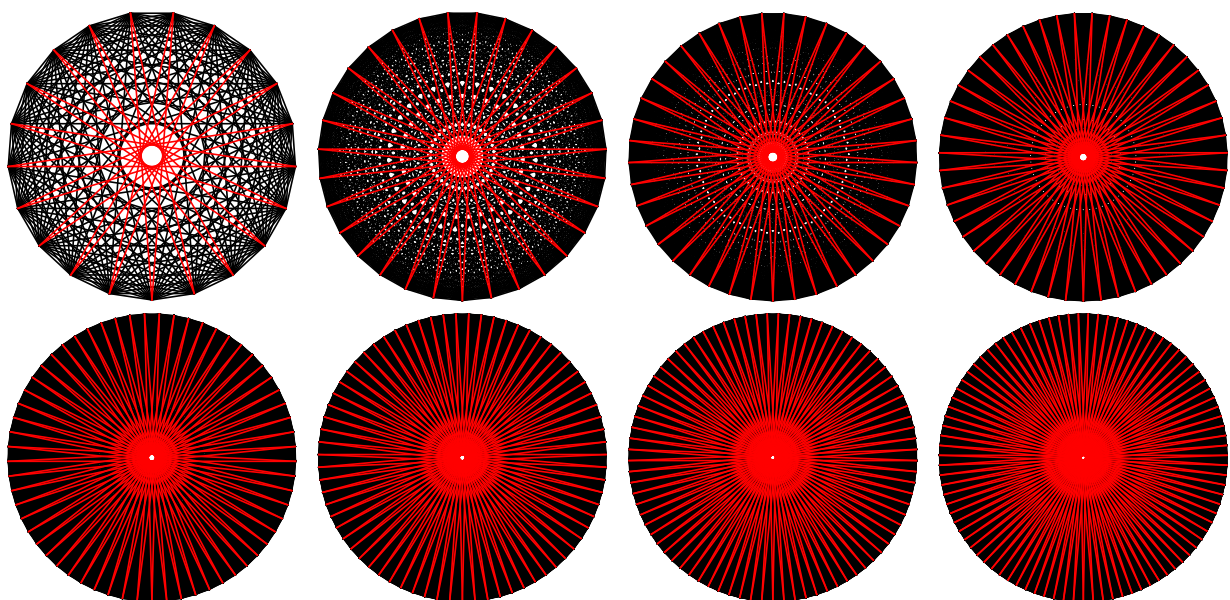

Figure 4. Largest small polygons $LSP(n)$ for $n = 21, 31, 41, 51, 61, 71, 81, 91$.

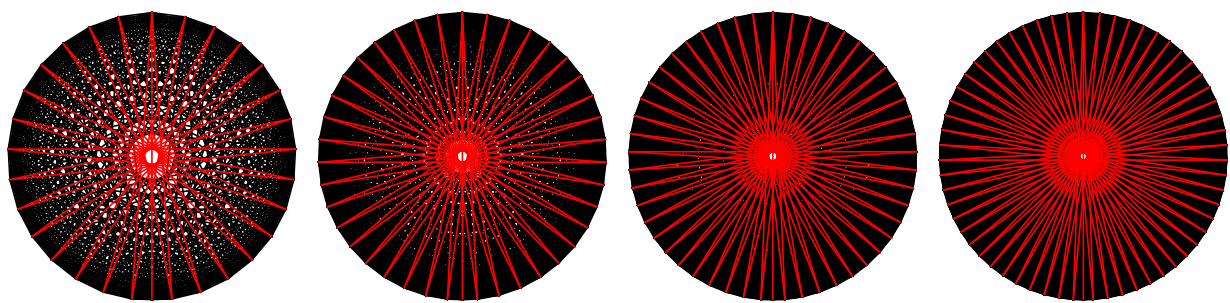

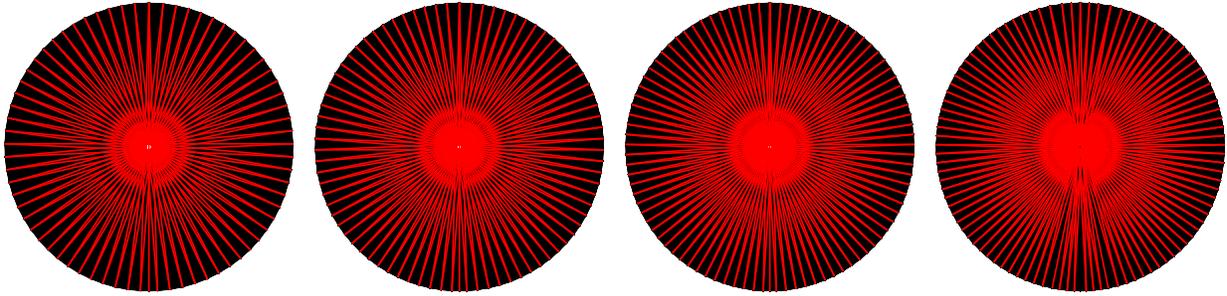

Figure 5. Largest small polygons $LSP(n)$ for $n = 30, 40, 50, 60, 70, 80, 90, 100$.

The next two figures summarize the difference between the area $A(n)$ of the optimized polygon and $\pi/4$ on as a function of $n$, on a loglog-plot scale.

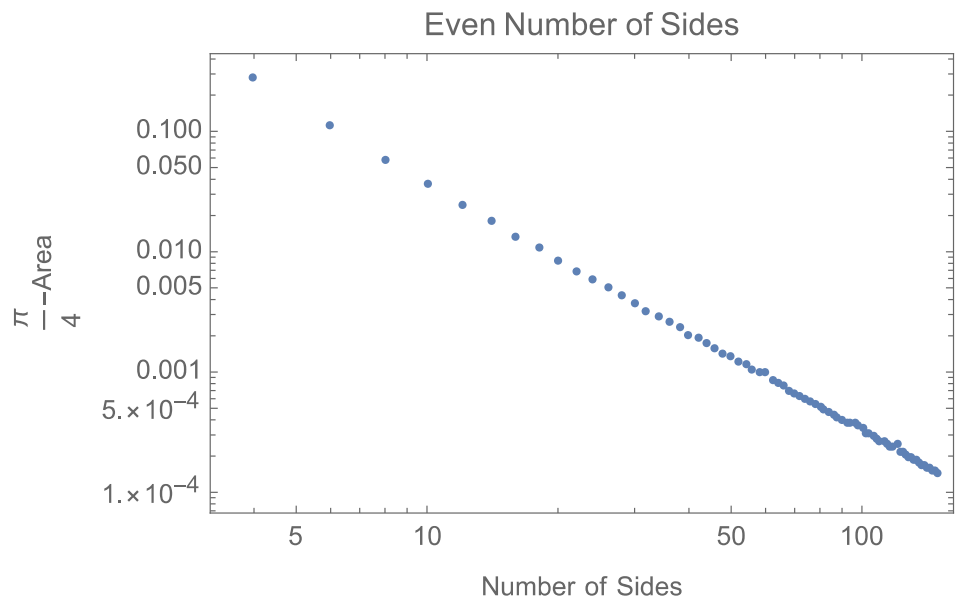

Figure 6. Difference between the area $A(n)$ of the optimized polygon and $\pi/4$, for selected even values $4 \leq n \leq 150$.

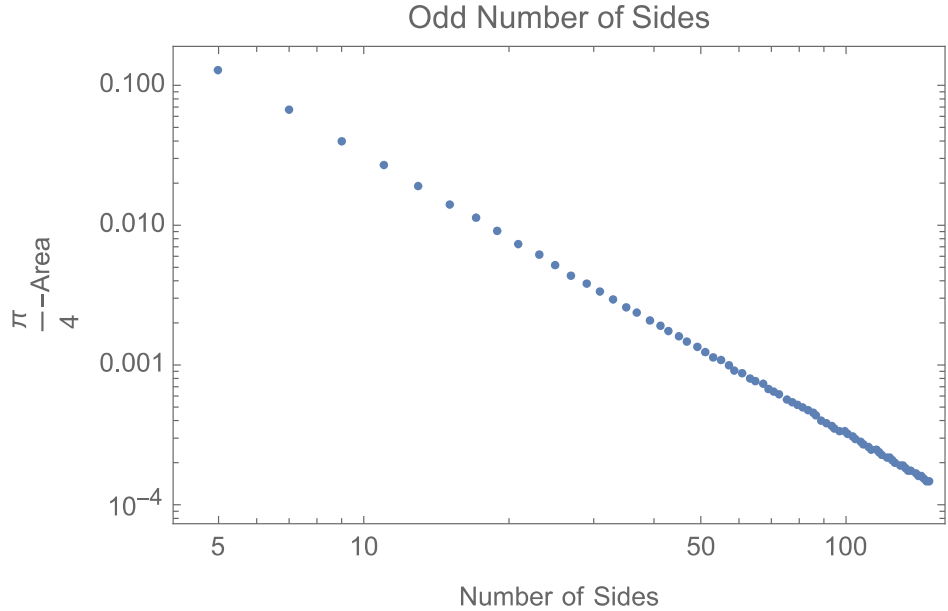

Figure 7. Difference between the area $A(n)$ of the optimized polygon and $\pi/4$, for selected odd values $5 \leq n \leq 149$.

Our calculations also reveal a small, but non-negligible difference between these figures: the numerically estimated slope of the plot for an even number of sides is approximately $-2.04618$, and the slope of the plot for an odd number of sides is approximately $-1.99848$.

We conclude the presentation of numerical results by emphasizing that the tightened LSP model offers advantages over the standard model. Specifically, IPOPT performs well on the tightened model, but it exhibits inferior performance on the standard model already for values $6 \leq n \leq 80$, as observed by Pintér (2020). Figure 8 illustrates superior IPOPT performance on tightened LSP model-instances up to $n = 300$, when compared to Figure 9 (standard LSP model with the same fixed starting solution as used in the tightened model) and Figure 10 (standard LSP model with $n$ random starting solutions). In the latter case, the solver runtimes also become longer: therefore, we conducted experiments only up to $n = 100$.

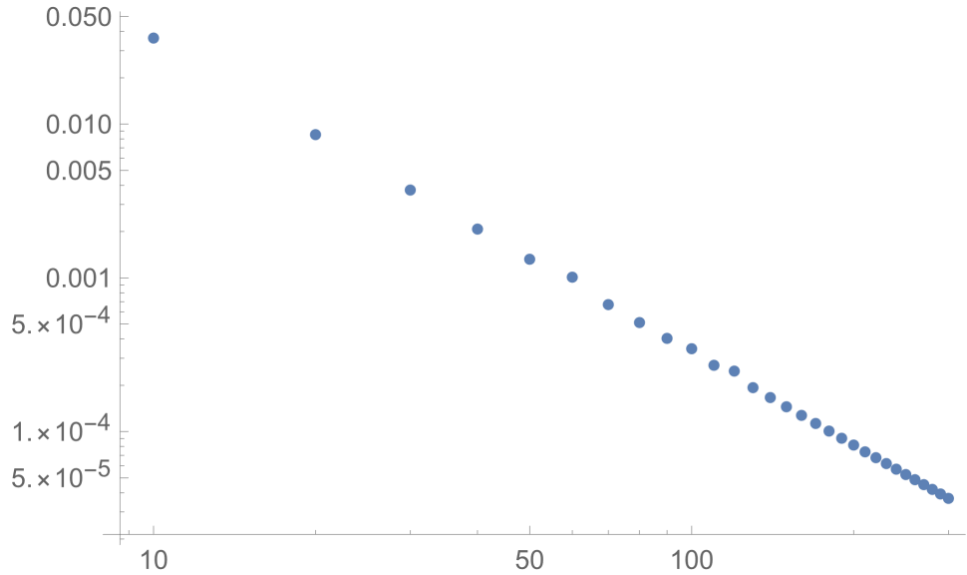

Figure 8. IPOPT performance on tightened LSP model-instances up to $n = 300$.
Difference between the area $A(n)$ of the optimized polygon and $\pi/4$.

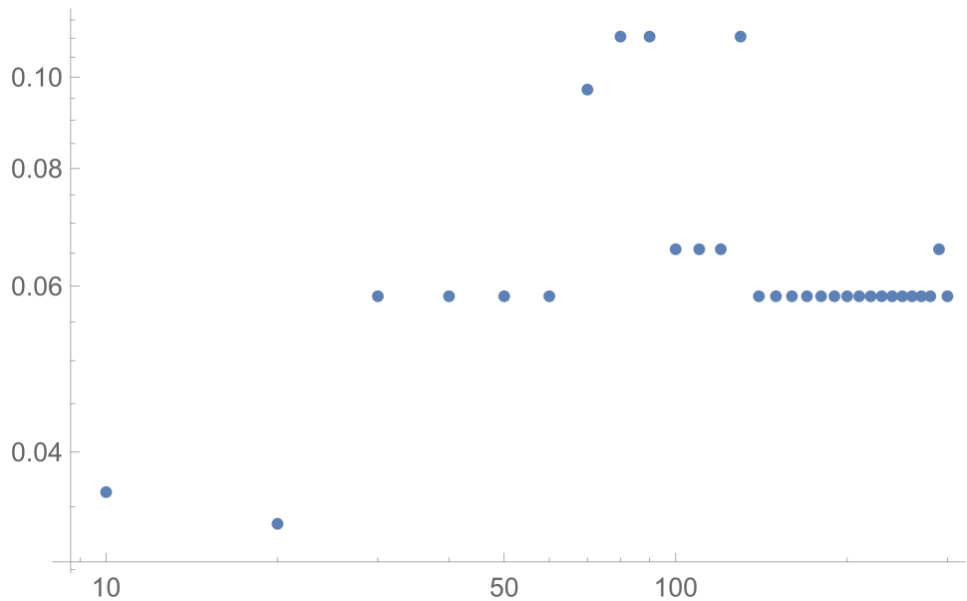

Figure 9. IPOPT performance on standard LSP model-instances, with the proposed starting solution, up to $n = 300$. Difference between the area $A(n)$ of the optimized polygon and $\pi/4$.

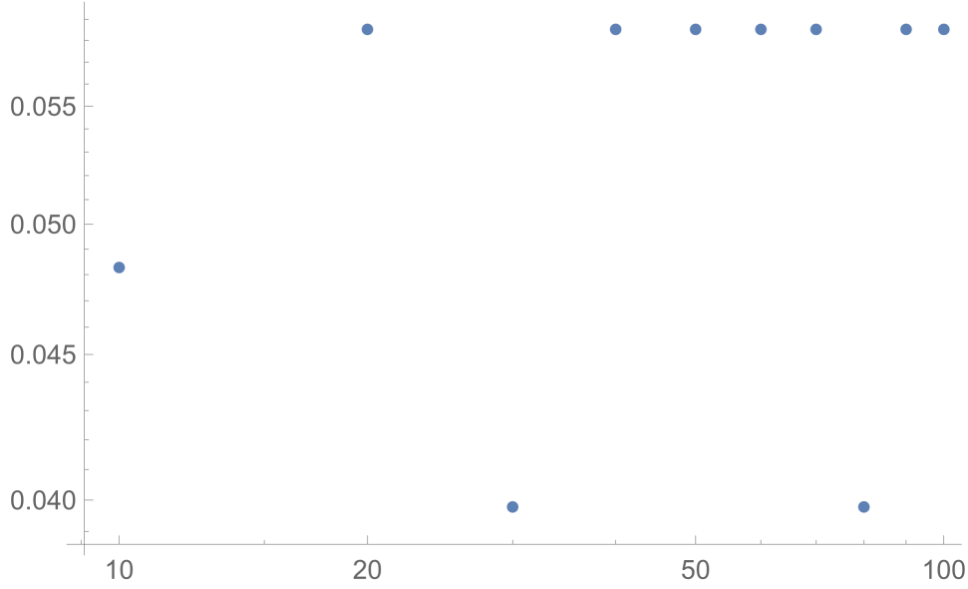
Figure 10. IPOPT performance on standard LSP model-instances, with random starting solutions up to $n=100$. Difference between the area $A(n)$ of the optimized polygon and $\pi/4$.

## 4 $LSP(n)$ Regression Models for Even and Odd $n$

As it is known, and also well illustrated by Figures 2 to 5, for large $n$ the optimal $LSP(n)$ configuration approaches the circle with unit diameter; hence the corresponding area limit is $A(\infty) = \pi/4 \sim 0.7853981634$. Comparing this limit value to our optimum estimate obtained for $A(1000) \sim 0.7853949284$, the ratio $A(1000)/(\pi/4)$ approximately equals $0.9999958811$. Hence, our $A(1000)$ estimate already leads to a fairly close approximation of the limit value.

Based on this observation and using the numerical results obtained, we developed the following regression models, separately for even and for odd values of $n$.

For even values of $n$, the regression model derived is
$$A(n) \sim \frac{\pi}{4} - 0.023182\left(\frac{1}{n}\right) - 2.630729\left(\frac{1}{n^2}\right) - 7.360373\left(\frac{1}{n^3}\right). \tag{11}$$

For odd values of $n$, the regression model derived is
$$A(n) \sim \frac{\pi}{4} - 0.024249\left(\frac{1}{n}\right) - 3.054853\left(\frac{1}{n^2}\right) - 0.131257\left(\frac{1}{n^3}\right). \tag{12}$$

Notice that the coefficients of the leading terms $(1/n)$ are fairly close to each other, with large differences between the other (paired) coefficients of $(1/n^2)$ and $(1/n^3)$. We received $p$-values (observed significance levels) well below $0.0001$ for all coefficients. This finding indicates that we have very strong statistical evidence suggesting that the regression coefficients are different from zero. Figures 11 and 12 depict the predicted results using models (11) and (12).

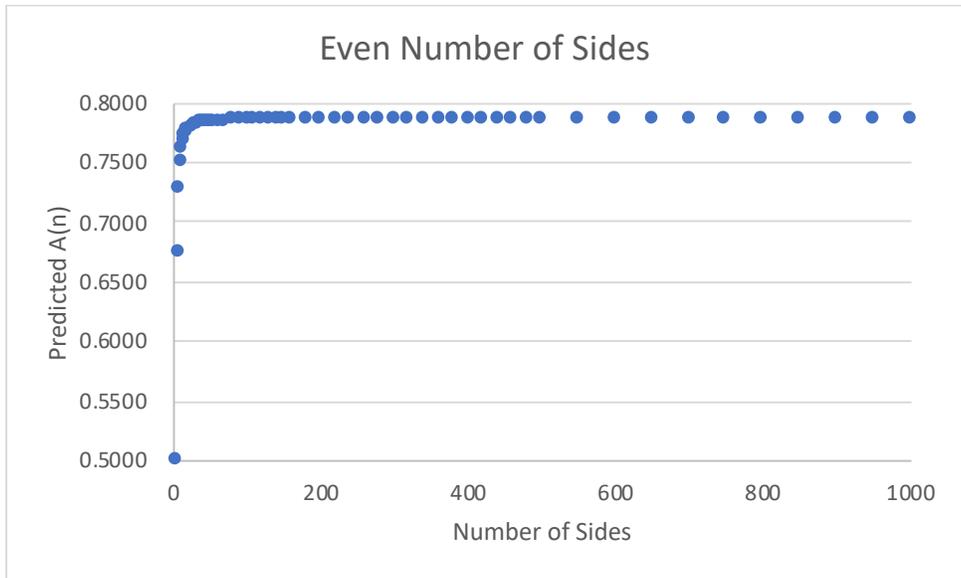

Figure 11. The nonlinear regression model (11) predicting $A(n)$ shown for $4 \leq n \leq 1000$.

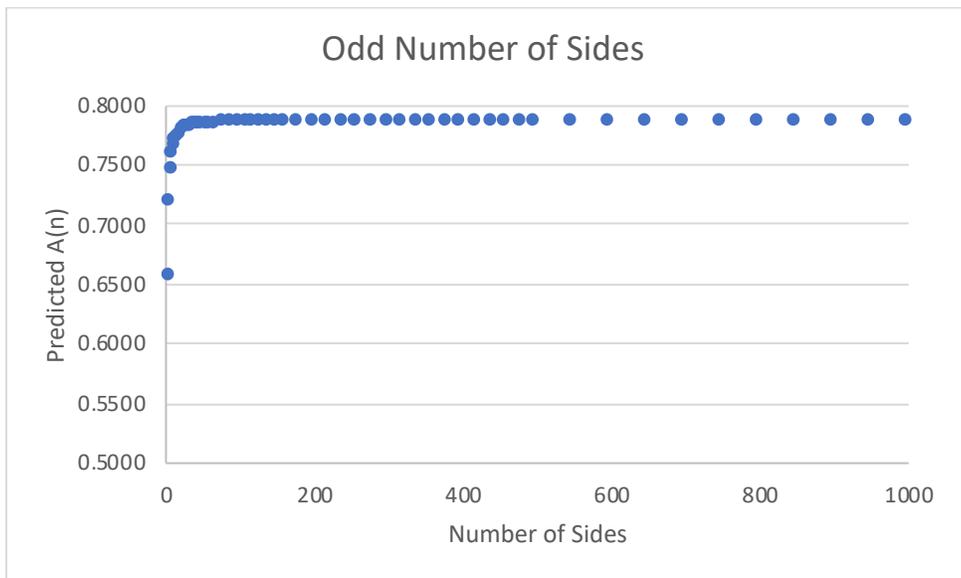

Figure 12. The nonlinear regression model (12) predicting $A(n)$ shown for $3 \leq n \leq 999$.

Since the two series closely overlap, the preceding two figures depicting observed *vs.* predicted $A(n)$ results become more useful when we zoom in and reduce the ranges considerably. In Figure 13, we display observed *vs.* predicted $A(n)$ results for $24 \leq n \leq 100$ and $0.779 \leq A(n) \leq 0.785$. The observed *vs.* predicted series appear visually different from each other only with high levels of magnification.

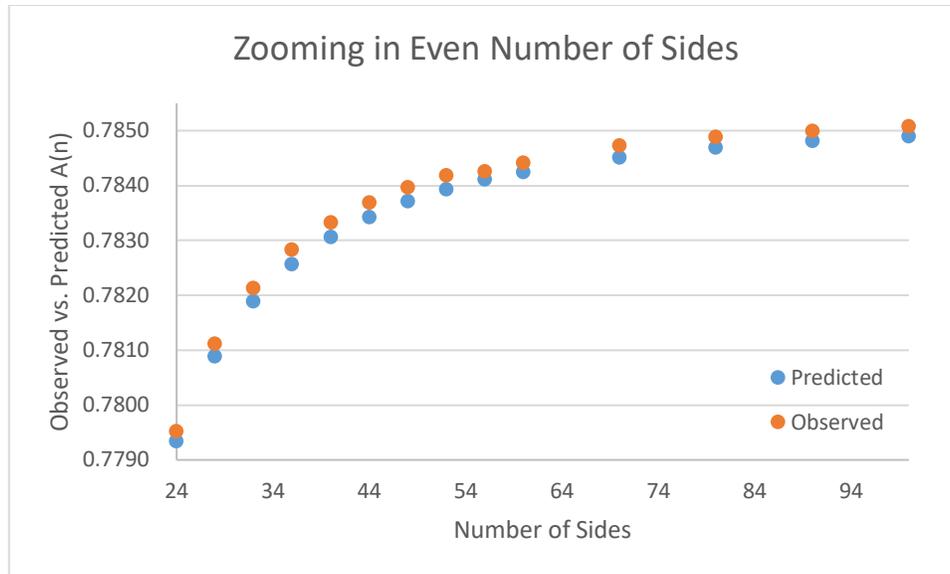

Figure 13. Zooming in: observed *vs.* predicted $A(n)$ results using regression model (11) for even $24 \leq n \leq 100$.

The regression models found can be used to directly estimate selected values from the entire sequence $\{A(n)\}$, including larger values of $n$ which have not been studied earlier and may be out of the range of current optimization solver capabilities. For example,

$$A(2000) \sim 0.7853872310, \qquad A(10000) \sim 0.7853958715.$$

$$A(1999) \sim 0.78538526834, \qquad A(9999) \sim 0.7853957077.$$

It is instructive to compare these estimates to $A(\infty) = \pi/4 \sim 0.7853981634$. Earlier numerical examples, with a different regression model based on results for even $6 \leq n \leq 80$, are presented in Pintér (2020).

## 5 Conclusions

Our study addresses the problem of finding numerically the sequence of largest small $n$-polygons $LSP(n)$ with unit diameter and maximal area $A(n)$, in principle aiming for all even (and odd) values of $n > 2$. This long-standing mathematical "puzzle" leads to an interesting class of nonlinear (global) optimization problems. First, we propose a revised LSP model, and later on we demonstrate its advantages compared to the standard model. Using the *Mathematica* modeling environment with the IPOPT solver option, and a new initial solution estimate, we can efficiently derive numerical solutions for a range of both even and odd values of $n$, up to $n \leq 1000$. Our results compare well to all the best results reported earlier for significantly lower values of $n$. We also propose regression models that enable the direct estimation of the optimal area sequence $\{A(n)\}$, for arbitrary even and odd values of $n$.

The LSP problem-class is one of those nice "puzzles" that can be described in a few words, yet lead to hard theoretical and numerical challenges. Therefore, this model-class, similarly to many other scientifically important configuration design models can be used also in software benchmarking tests. We think that such problems can serve as a significant addendum to the collection of (well-frequented, and often much simpler) unconstrained or box-constrained test problems. For further scalable and increasingly hard-to-solve object configuration models, we refer to some of our earlier and more recent studies: Castillo, Kampas, and Pintér (2008), Kampas, Castillo, and Pintér (2019), Kampas, Pintér, and Castillo (2020), Pintér (2001), Pintér, Kampas, and Castillo (2017).

# References

Castillo, I., Kampas, F.J. and Pintér, J.D. (2008) Solving circle packing problems by global optimization: numerical results and industrial applications. *European Journal of Operational Research* 191 (2008) 786–802.

Mossinghoff, M.J. (2006), Isodiametric problems for polygons. *Discrete and Computational Geometry*, 36 (2), pp. 363–379.

Kampas, F.J., Castillo, I., and Pintér, J.D. (2019), Optimized ellipse packings in regular polygons. *Optimization Letters* 13, 1583–1613.

Kampas, F.J., Pintér, J.D. and Castillo, I. (2020), Packing ovals in optimized regular polygons. *Journal of Global Optimization* 77, 175–196.

Pintér, J.D. (2001), Globally optimized spherical point arrangements: model variants and illustrative results. *Annals of Operations Research* 104, 213-230.

Pintér, J.D. (2020), Largest small n-polygons: numerical optimum estimates for n ≥ 6. To appear in: Al-Baali, M., Grandinetti, L., Purnama, A., Eds. *Recent Developments in Numerical Analysis and Optimization NAO-V*, Muscat, Oman, January 2020. Springer Nature, New York, 2021.

Pintér, J.D. and Kampas, F.J. (2005), O.R. model development and optimization with *Mathematica*. In: Golden, B., Raghavan, S., and Wasil, E., Eds. *The Next Wave in Computing, Optimization, and Decision Technologies*, pp. 285-302. Springer Science + Business Media, New York.

Pintér, J.D. and Kampas, F.J. (2013), Benchmarking nonlinear optimization software in technical computing environments: global optimization in *Mathematica* with *MathOptimizer Professional*. *TOP* 21, 133-162.

Pintér, J.D., Kampas, F.J. and Castillo, I. (2017), Globally optimized packings of non-uniform size spheres in $R^d$: a computational study. *Optimization Letters* 12, 585–613.

Reinhardt, K. (1922), Extremale polygone gegebenen durchmessers, *Jahresbericht der Deutschen Mathematiker-Vereinigung* 31, 251–270.

Wächter, A. and Laird, C. (2020), IPOPT (Interior Point Optimizer), an open source software package for large-scale nonlinear optimization. https://coin-or.github.io/IPOPT/. The contributors to the IPOPT Project are listed at https://coin-or.github.io/IPOPT/AUTHORS.html. (Retrieved on August 29, 2020.)

Weisstein, E.W. (2020), *"Biggest Little Polygon."* From *MathWorld* – A Wolfram Web Resource. http://mathworld.wolfram.com/BiggestLittlePolygon.html. (Retrieved on August 29, 2020.)

Wolfram Research (2020a), *Mathematica*, Current release 12.1 (as of August 2020). Wolfram Research, Champaign, IL.

Wolfram Research (2020b), Optimizing with IPOPT. https://reference.wolfram.com/language/IPOPTLink/tutorial/OptimizingWithIPOPT.html.  (Retrieved on August 29, 2020.)